\documentclass[12pt,a4paper]{article}
\usepackage{amssymb,amsmath}
\usepackage{amsfonts}
\usepackage{epsfig}
\usepackage{float}
\usepackage{graphicx}
\usepackage{color}
\usepackage[lined,boxed,commentsnumbered]{algorithm2e}
\usepackage{multicol}
\newcommand{\para}{\par\vspace{.25cm}}

\newtheorem{theorem}{Theorem}
\newtheorem{lemma}{Lemma}
\newtheorem{cor}{Corollary}

\begin{document}
\baselineskip 18pt \title{\bf Extremely strong Shoda pairs with \texttt{GAP}}
\author{ Gurmeet K. Bakshi and Sugandha Maheshwary{\footnote {Research supported by CSIR, India (File No: SPM-09/135(0107)/2011-EMR-I) is gratefully acknowledged} \footnote{Corresponding author}} \\ {\em \small Centre for Advanced Study in
Mathematics,}\\
{\em \small Panjab University, Chandigarh 160014, India}\\{\em
\small email: gkbakshi@pu.ac.in and msugandha.87@gmail.com } }
\date{}
{\maketitle}
\begin{abstract}\noindent
{We provide algorithms to compute a complete irredundant set of extremely strong Shoda pairs of a finite group $G$ and the set of primitive central idempotents of the rational group algebra $\mathbb{Q}[G]$ realized by them. These algorithms are also extended to write new algorithms for computing a complete irredundant set of strong Shoda pairs of $G$ and the set of primitive central idempotents of $\mathbb{Q}[G]$ realized by them. Another algorithm to check whether a finite group $G$ is normally monomial or not is also described.}
\end{abstract}\vspace{.25cm}
 \vspace{.25cm} 

\section{Introduction}

Let $G$ be a finite group and let $\mathbb{Q}[G]$ be the rational group algebra of $G$. A \emph{strong Shoda pair} of $G$, introduced by Olivieri, del R{\'{\i}}o and Sim{\'o}n \cite{Oli}, is a pair $(H,K)$ of subgroups of $G$ with the subgroups $H$ and $K$ satisfying some technical conditions. In \cite{BM}, a strong Shoda pair $(H,K)$ with $H$ normal in $G$ is called as an \emph{extremely strong Shoda pair} of $G$. An important property (\cite{Oli}, Proposition 3.3) of the strong Shoda pairs of $G$ is that each such pair $(H,K)$ determines a \emph{primitive central idempotent} of $\mathbb{Q}[G]$, called the \emph{primitive central idempotent of} $\mathbb{Q}[G]$  \emph{realized by} $(H,K)$, which is denoted by $e(G,H,K)$. Let $E$ be the set of all primitive central idempotents of $\mathbb{Q}[G]$ and $E_{SSP}$ (resp. $E_{ESSP}$) be the set of primitive central idempotents of $\mathbb{Q}[G]$ realized by the strong Shoda pairs (resp. extremely strong Shoda pairs) of $G$. The groups $G$ for which $E=E_{SSP}$ are called \emph{strongly monomial} groups and are known to constitute a large class of monomial groups, including abelian-by-supersolvable groups \cite{Oli}. Also, in \cite{BM}, it has been proved that $E=E_{SSP}=E_{ESSP}$ if, and only if, $G$ is a \emph{normally monomial} group i.e., every complex irreducible character of $G$ is induced from a linear character of a normal subgroup of $G$. The \texttt{GAP}~\cite{GAP} package \texttt{Wedderga} \cite{Wedd} features the function \texttt{PrimitiveCentralIdempotentsByStrongSP(QG);} that computes the set $E_{SSP}$ for the rational group algebra  $\mathbb{Q}[G]$ and the function \texttt{StrongShodaPairs(G);} that determines a subset $X$ of strong Shoda pairs of $G$ such that $(H,K)\mapsto e(G,H,K)$ defines a bijection from $X$ to $E_{SSP}$. Such a set $X$ is called a \emph{complete irredundant set of strong Shoda pairs} of $G$. These functions are based on the search algorithms provided by Olivieri and del R{\'{\i}}o \cite{OliA}. Another relevant feature of \texttt{Wedderga} is the function \texttt{IsStronglyMonomial(G);} which checks whether the group $G$ is strongly monomial or not. Using this function, it has been revealed in \cite{Olte} that all the monomial groups of order less than $1000$ are strongly monomial.\\

In this paper, we provide an algorithm to compute a complete irredundant set of extremely strong Shoda pairs of $G$. This algorithm is based on the work in~\cite{BM}. We further extend this algorithm by combining it with the search algorithm provided by Olivieri and del R{\'{\i}}o \cite{OliA} to obtain a new algorithm that computes a complete irredundant set of strong Shoda pairs of $G$. As a consequence, we obtain algorithms to write the sets $E_{ESSP}$ and $E_{SSP}$ of primitive central idempotents of $\mathbb{Q}[G]$ realized by extremely strong Shoda pairs of $G$ and those re alized by strong Shoda pairs of $G$ respectively. Another algorithm to check whether a finite group~$G$ is normally monomial or not also follows as a consequence. These algorithms are given in Section 3 and enable us to write the following functions in \texttt{GAP} language:

 \begin{itemize}
  \item \texttt{ExtStrongShodaPairs(G);} which computes a complete irredundant set of extremely strong Shoda pairs of $G$ i.e., a subset $X$ of extremely strong Shoda pairs of $G$, such that such that $(H,K)\mapsto e(G,H,K)$ gives a bijection from $X$ to $E_{ESSP}$.
      \item \texttt{StShodaPairs(G);} which computes a complete irredundant set of strong Shoda pairs of $G$.
  \item \texttt{PrimitiveCentralIdempotentsByExtSSP(QG);} which computes the set of primitive central idempotents of $\mathbb{Q}[G]$ realized by extremely strong Shoda pairs of $G$.
  \item \texttt{PrimitiveCentralIdempotentsByStSP(QG);} which computes the set of primitive central idempotents of $\mathbb{Q}[G]$ realized by strong Shoda pairs of $G$.
  \item \texttt{IsNormallyMonomial(G);} which checks whether the group $G$ is normally monomial or not.
\end{itemize}

 Using the function \texttt{IsNormallyMonomial(G);} we have searched for normally monomial groups among the groups in \texttt{GAP} library of small groups. The search indicates that the class of normally monomial groups is a substantial class of monomial groups. It may also be mentioned that if $G$ is a normally monomial group, then the output obtained by the functions \texttt{StShodaPairs(G);} and \texttt{PrimitiveCentralIdempotentsByStSP(QG);} is same as that obtained by \linebreak \texttt{ExtStrongShodaPairs(G);} and \texttt{PrimitiveCentralIdempotentsByExtSSP(QG);} respectively. Furthermore, for a finite group $G$, the functions \texttt{StShodaPairs(G);} and \texttt{PrimitiveCentralIdempotentsByStSP(QG);}  are alternative to the functions \texttt{StrongShodaPairs(G);} and \texttt{PrimitiveCentralIdempotentsByStrongSP(QG);} respectively, which are currently available in \texttt{Wedderga}.

 \para In Section 4, we compare the runtimes of the function \texttt{StShodaPairs(G);} with \texttt{StrongShodaPairs(G);} for a large and evenly spread sample of groups of order up to $2000$. For this sample, the functions \texttt{PrimitiveCentralIdempotentsByStSP(QG);} and \texttt{PrimitiveCentralIdempotentsByStrongSP(QG)}; are also compared for runtimes. It is observed that these new functions show significant improvement in the time taken to compute the same outputs. Further, in order to observe the performance separately for solvable and non solvable groups, we also compared the runtimes of \texttt{StShodaPairs(G);} with \texttt{StrongShodaPairs(G);} for another two samples. The sample of solvable groups consists of all the groups of odd order up to 2000, and the other sample consists of all non solvable groups of order up to 2000. It is observed that the performance of \texttt{StShodaPairs(G);} is exceptionally better in comparison with that of \texttt{StrongShodaPairs(G);} for solvable groups. However, in the case of non solvable groups, the performance of the two functions is almost identical. Finally, we describe the reasons for the difference in the performance of these functions.

\section{Notation and Preliminaries}
Throughout this paper, $G$ denotes a finite group. By $H \leq G$ (resp. $H\unlhd G$), we mean that $H$ is a subgroup (resp. normal subgroup) of $G$. For $H \leq G$, $[G:H]$ denotes the index of $H$ in $G$, $N_{G}(H)$ denotes the normalizer of $H$ in $G$, $\operatorname{core}_{G}(H)=\bigcap_{x\in G}xHx^{-1}$ and $\hat{H}=\frac{1}{|H|}\sum_{h \in H}h$, where $|H|$ is the order of $H$. For $K\unlhd H \leq G$, write
 $$\varepsilon (H,K) := \begin{cases} \hat{H}, & {\rm if}~ H = K; \\ \prod (\hat{K}-\hat{L}), &  {\rm otherwise,} \end{cases}$$ where $L$ runs over the minimal normal subgroups of $H$ containing $K$ properly. Set \begin{center}$e(G,H,K)$ := the sum of all the distinct $G$-conjugates of $\varepsilon(H,K)$. \end{center}
 Let $\varphi$ denote the Euler phi function.
Denote by $\operatorname{Irr}(G)$, the set of all complex irreducible characters of $G$. For $\chi \in \operatorname{Irr}(G)$, $\mathbb{Q}(\chi)$ denotes the field obtained by adjoining to $\mathbb{Q}$, all the character values $\chi(g),~g \in G$, and $\operatorname{Gal}(\mathbb{Q}(\chi)/\mathbb{Q})$ is the Galois group of the extension $\mathbb{Q}(\chi)$ over $\mathbb{Q}$.
\para It is well known that $\chi \mapsto e_{\mathbb{Q}}(\chi):=\frac{\chi(1)}{|G|}\sum_{\sigma \in \operatorname{Gal}(\mathbb{Q}(\chi)/\mathbb{Q}) }\sum_{g \in G}\sigma(\chi(g^{-1}))g$ defines a surjective map from $\operatorname{Irr}(G)$  to the set of primitive central idempotents of the rational group algebra $\mathbb{Q}[G]$. If $\chi$ is the trivial character of $G$, then it is easy to see that $e_{\mathbb{Q}}(\chi)=\hat{G}$.\\

Olivieri et al \cite{Oli} proved the following:
\begin{theorem}\label{t1}{\rm({\cite{Oli}, Lemma 1.2, Theorem 2.1})}
\begin{enumerate}
  \item If $\chi$ is a non trivial linear character of $G$ with kernel $N$ then
   \begin{equation}\label{e1}
    e_{\mathbb{Q}}(\chi)=\varepsilon(G,N).
 \end{equation}
  \item If $\chi$ is monomial, i.e. $\chi$ is induced from a linear character $\psi$ of a subgroup $H$ of $G$, then there exists $\alpha \in  \mathbb{Q}$ such that
      \begin{equation}
        e_{\mathbb{Q}}(\chi)= \alpha e(G, H, K),      \end{equation}
       where $K$ is the kernel of the character $\psi$. Furthermore $\alpha=1$, if the distinct $G$-conjugates of $\varepsilon(H,K)$ are mutually orthogonal.

      \end{enumerate}
\end{theorem}
\para Shoda (see \cite{curt}, Theorem 45.2) gave a criteria to decide if the character of $G$ induced by a linear character $\psi$ of a subgroup $H$ of $G$ is irreducible, in terms of $H$ and the kernel $K$ of $\psi$. A pair $(H,K)$ of $G$ satisfying
this criteria is called a \emph{Shoda pair} (\cite{Oli}, Definition 1.4) of $G$. A \emph{strong Shoda pair} (\cite{Oli}, Definition 3.1) of $G$ is a pair $(H,K)$ of subgroups satisfying the following conditions:
\begin{quote}\begin{description}
  \item[(i)] $K \unlhd H \unlhd N_{G}(K)$;
  \item[(ii)] $H/K$ is cyclic and a maximal abelian subgroup of $N_{G}(K)/K$;
  \item[(iii)] the distinct $G$-conjugates of $\varepsilon(H,K)$ are mutually orthogonal.
\end{description}\end{quote}
As the name suggests, each strong
Shoda pair of $G$ is also a Shoda pair of $G$ (\cite{Oli}, Proposition 3.3).
A strong Shoda pair $(H,K)$ of $G$ is called an \emph{extremely strong Shoda pair} of $G$, if $H\unlhd G$. Observe that $(G,G)$ is always an extremely strong Shoda pair of $G$.
\para From Theorem \ref{t1}, it follows that if $(H, K)$ is a strong Shoda pair of $G$, then $e(G,H,K)$ is a primitive central idempotent of $\mathbb{Q}[G]$, called {\it the primitive central idempotent of $\mathbb{Q}[G]$ \it realized by $(H,K)$}. For a strong Shoda pair $(H,K)$ of $G$, we denote by $\operatorname{dim}(H,K)$, the $\mathbb{Q}$-dimension of the simple component $\mathbb{Q}[G]e(G,H,K)$ of $\mathbb{Q}[G]$. In view of (\cite{Oli}, Proposition 3.4), $\operatorname{dim}(H,K)$ equals $\varphi([H:K])[N_{G}(K):H][G:N_{G}(K)]^{2}.$ Two strong (resp. extremely strong) Shoda pairs $(H_{1}, K_{1})$ and $(H_{2}, K_{2})$ of $G$ are said to be {\it equivalent} if  $e(G, H_{1}, K_{1})=e(G, H_{2}, K_{2})$. A complete set of representatives of distinct equivalence classes of strong  (resp. extremely strong) Shoda pairs of $G$ is called a {\it complete irredundant set of strong {\rm (}resp. extremely strong{\rm)} Shoda pairs} of $G$.
\para We now recall the method given in \cite{BM} to compute a complete irredundant set of extremely strong Shoda pairs of a finite group $G$.\para Let $\mathcal{N}$ be the set of all the distinct normal subgroups of $G$. For $N \in \mathcal{N}$, let $A_{N}$ be a normal subgroup of $G$ containing $N$ such that $A_{N}/N$ is an abelian normal subgroup of maximal order in $G/N$. Note that the choice of $A_{N}$ is not unique. However, we need to fix one such $A_{N}$. For a fixed $A_{N}$, set \vspace{.5cm}\\ $\begin{array}{lll}
  \mathcal{D}_{N}: & {\rm the ~set ~of ~all ~subgroups~} D ~{\rm  of} ~A_{N}~ {\rm containing ~} N~{\rm  such ~that~}
   \operatorname{core}_{G}(D)=N, \\&  A_{N}/D ~{\rm    is ~ cyclic ~  and ~  is ~a ~ maximal~ abelian ~ subgroup~  of~}   N_{G}(D)/D.  \vspace{.2cm}\\ \mathcal{T}_{N}: &  {\rm a ~set~ of ~representatives~ of~ } \mathcal{D}_{N} {\rm ~  under ~the~ equivalence~ relation~ defined~ by} \\ & {\rm
    conjugacy~ of~ subgroups~in~} G. \vspace{.2cm} \\  \mathcal{S}_{N}: & \{( A_{N},D)~|~ D \in \mathcal{T}_{N}\}.\end{array}$ \vspace{.2cm}

\begin{theorem}\label{t2}{\rm(\cite{BM}, Theorem 1)} Let $G$ be a finite group. Then, \\ {\rm(i)}   $\cup_{N \in \mathcal{N} }\mathcal{S}_{N} $ is a complete irredundant set of extremely strong Shoda pairs of $G$. \\{\rm(ii)}  $\{ e(G, A_{N}, D)\, |\, ( A_{N},D)  \in \mathcal{S}_{N},~ N \in \mathcal{N} \}$ is a complete set of primitive central idempotents of $\mathbb{Q}[G]$ if, and only if, $G$ is normally monomial.
 \end{theorem}

 \para It may be noted that in Theorem \ref{t2}, the choice of $A_{N}$ is irrelevant. For $N \in \mathcal{N}$, let $A_{N}'$ be another normal subgroup of $G$ containing $N$ such that $A_{N}'/N$ is an abelian normal subgroup of maximal order in $G/N$ and let $\mathcal{D}_{N}'$, $\mathcal{T}_{N}'$ and $\mathcal{S}_{N}'$ be defined corresponding to $A_{N}'$. Then any pair in $\mathcal{S}_{N}'$ is equivalent to a pair in $\mathcal{S}_{N}$ and vice versa. This is because, if $(A_{N}',D')\in \mathcal{S}_{N}'$ and $\psi$ is a linear character of $A_{N}'$ with kernel $D'$, then $\psi^{G}$ is irreducible and hence by (\cite{BM}, Lemma 1), there exists $(A_{N},D)\in \mathcal{S}_{N}$ such that $e_{\mathbb{Q}}(\psi^{G})=e(G,A_{N},D).$ However, in view of Theorem \ref{t1}, $e_{\mathbb{Q}}(\psi^{G})=e(G,A_{N}',D').$ This gives that $(A_{N},D)$ is equivalent to $(A_{N}',D')$. The reverse conclusion holds similarly.

\begin{cor}{\rm(\cite{BM}, Corollary 1)}\label{c1}  If $G$ is a normally monomial group, then $\bigcup_{N \in \mathcal{N}}\mathcal{S}_{N}$ is  a complete irredundant set of strong Shoda pairs of $G$.
\end{cor}

\begin{cor}\label{c2}{\rm(\cite{BM}, Corollary 2)} A finite group $G$ is normally monomial if, and only if,
  $$  \sum_{N\in \mathcal{N}}\sum_{(A_{N},D) \in \mathcal{S}_{N}}dim(A_{N},D) = |G|.$$

\end{cor}

\section{Algorithms}
\par We shall use the notation developed in the previous section.

 \subsection{Extremely Strong Shoda Pairs}
We provide Algorithm 1, which computes the set $ESSP$, which is a complete irredundant set of extremely strong Shoda pairs of a given finite group $G$. This algorithm is based on Theorem \ref{t2}. It mainly requires the set $\mathcal{N}$ of normal subgroups of $G$ and the computation of $\mathcal{S}_{N}$ for each $N\in \mathcal{N}$. The set $\mathcal{S}_{N}$ is computed as explained in Section 2 and by using Lemmas \ref{l1}-\ref{l2} to avoid unnecessary computations.
\begin{lemma}\label{l1} For a normal subgroup $N$ of $G$, the following hold:
\begin{description}
  \item[(i)] If $G/N$ is abelian, then
  $$\mathcal{S}_{N} = \begin{cases} \{(G, N)\}, ~& {\rm if}~G/N~{\rm is~  cyclic;}\\ \emptyset, & {\rm otherwise.}
    \end{cases}$$
  \item[(ii)] If $G/N$ is non abelian and $A_{N}/N$ is cyclic, then
  $$\mathcal{S}_{N} = \begin{cases} \{(A_{N}, N)\}, & ~{\rm if}~ A_{N}/N ~{\rm is~ a~ maximal~ abelian~ subgroup ~of} ~G/N;\\ \emptyset, & {\rm otherwise.}
    \end{cases}$$
\end{description}

\end{lemma}
{\bf Proof.} Follows immediately from the definition of $\mathcal{S}_{N}$. $\Box$

\begin{lemma}\label{l3} If $\mathcal{M} \subseteq \mathcal{N}$ is such that $$ \sum_{N\in \mathcal{M}}\sum_{(A_{N},D) \in \mathcal{S}_{N}}dim(A_{N},D)=|G|,$$
then $\mathcal{S}_{N}= \emptyset~ $ for all $~N \in \mathcal{N}\setminus \mathcal{M}$.\end{lemma}
{\bf Proof.} The primitive central idempotents $e(G,A_{N},D)$ for $(A_{N},D)\in \mathcal{S}_{N} $, $N \in \mathcal{N}$, are distinct, as $\bigcup_{N \in \mathcal{N}}\mathcal{S}_{N}$ is a complete irredundant set of extremely strong Shoda pairs of $G$. Therefore, $\bigoplus_{N \in \mathcal{N}}\bigoplus_{(A_{N},D)\in \mathcal{S}_{N}}\mathbb{Q}[G]e(G,A_{N},D)$
 is a direct summand of $\mathbb{Q}[G]$, and hence its $\mathbb{Q}$-dimension is at most $|G|$. Consequently,
\begin{eqnarray*}
  |G| &\geq & \sum_{N\in \mathcal{N}}\sum_{(A_{N},D)\in \mathcal{S}_{N}} dim(A_{N},D) \nonumber \\
   &\geq & \sum_{N\in \mathcal{M}}\sum_{(A_{N},D)\in \mathcal{S}_{N}} dim(A_{N},D) ~~~~(\because\mathcal{M} \subseteq \mathcal{N}) \nonumber\\
   &=& |G|.
\end{eqnarray*}
This yields that $\mathcal{S}_{N} =\emptyset$ for all $N\not \in \mathcal{M}$ and completes the proof.$~\Box$

\begin{lemma}\label{l2}
If $(H,K)$ is a strong Shoda pair of $G$ with $N=\operatorname{core}_{G}(K)$, then the centre of $G/N$ must be cyclic.\\
\indent In particular, if $N \in \mathcal{N}$  is such that  the centre of $G/N$ is not cyclic, then $\mathcal{S}_{N}= \emptyset$.
\end{lemma}
{\bf Proof.}  Let $aK$ be a generator of $H/K$ and let $\zeta$ be a primitive $m^{th}$ root of unity, where $m=[H:K]$. Consider the linear representation $\rho: H \rightarrow \mathbb{C}$ given by $x\mapsto \zeta^{i}$, if $xK=a^{i}K$, for $x \in H$. Since $(H,K)$ is a strong Shoda pair, $\rho^{G}$ is an irreducible representation of $G$. Now, as $\operatorname{ker}(\rho^{G})=\bigcap_{x\in G}x(\operatorname{ker}\rho)x^{-1}=\bigcap_{x\in G}xKx^{-1}=\operatorname{core}_{G}(K)=N$, the result follows from (\cite{IM}, Lemma 2.27).$~\Box$\\

\begin{algorithm}\label{A1}
\begin{scriptsize}
\KwData{A finite group $G$.}
$\mathcal{N}:=$Normal subgroups of $G$ (in decreasing order)\;
$ESSP:=[[G,G]]$\;
$SumDim:=1$\;
\While{$SumDim \neq |G|$}{$\mathcal{N}_{1}:=\emptyset$\;\For{$N$ in $\mathcal{N}$, $N\neq G$}{\eIf{$G'\subseteq N$}{\If{$G/N$ is cyclic}{Add $[G,N]$ to the list $ESSP$\;$dim:=dim(G,N)$; $SumDim:=SumDim+dim;$}}{Add $N$ to the list $\mathcal{N}_{1}$\;}}$\mathcal{N}_{2}:=\emptyset$\; \For{$N$ in $\mathcal{N}_{1}$}{\If{$Centre(G/N)$ is cyclic}{Add $N$ to the list $\mathcal{N}_{2}$\;}}$List0:=\emptyset$\; \For{$N$ in $\mathcal{N}_{2}$}{$A:=$ a normal subgroup of $G$ containing $N$ such that $A/N$ is an abelian normal subgroup of maximal order in $G/N$\;\If{$A\neq N$}{\eIf {$A/N$ is cyclic;  }{\If{$A/N$ is maximal abelian subgroup of $G/N$}{Add $[A,N]$ to the list $ESSP$\;$dim:=dim(A,N)$; $SumDim:=SumDim+dim;$}}{Add $[A,N]$ to the list $List0$\;}}}$LIST:=List0$\;
     \While{$LIST$ is non empty}{$A:=LIST[1][1]$\;$NA$:=Normal subgroups $D$ of $A$ such that $A/D$ is cyclic\;$LIST0$:=Pairs $p=[p[1],p[2]]$ in $LIST$ such that $p[1]=A$\;
     \For{$q=[q[1],q[2]]$ in $LIST0$}{$\mathcal{D}:=$Subgroups $D\in NA$ such that $\operatorname{core}_{G}(D)=q[2]$\;$\mathcal{T}:=$Distinct conjugate representatives of $\mathcal{D}$\;
     \For{$T$ in $\mathcal{T}$}{\If{$A/T$ is maximal abelian subgroup of $N_{G}(T)/T$}{  Add $[A,T]$ to the list $ESSP$\; $dim:=dim(A,T);$ $SumDim:=SumDim+dim;$}}
     }}

      $LIST:=LIST\setminus LIST0$;}

\KwResult{$ESSP$ }

\end{scriptsize}

\caption{Extremely strong Shoda pairs of $G$}

\end{algorithm}
We now describe Algorithm 1. The first step of the algorithm is to compute the list $\mathcal{N}$ of all the normal subgroups of $G$ in decreasing order. If $N=G$, we have $\mathcal{S}_{N}= \{(G,G)\}$. Therefore, we initially set the list $ESSP$, which is the list of extremely strong Shoda pairs of $G$ found at any stage of computation, to be \texttt{[[G,G]]}. As $\operatorname{dim}(G,G)=1$, we set $SumDim$, which denotes the sum of $\mathbb{Q}$-dimensions of the simple components of $\mathbb{Q}[G]$ corresponding to the elements in $ESSP$, to be 1. For $N\in \mathcal{N}$, $N \neq G$, if $N$ contains the commutator subgroup $G'$ of $G$, then the corresponding set $\mathcal{S}_{N}$ is computed using Lemma \ref{l1}(i). Otherwise, $\mathcal{S}_{N}$ is computed using Theorem \ref{t2} along with Lemmas \ref{l1} and \ref{l2}. In either of the two cases, if $\mathcal{S}_{N}\neq \emptyset$, then the elements of $\mathcal{S}_{N}$ are added to the list $ESSP$. Also, the sum of $\mathbb{Q}$-dimensions of simple components of $\mathbb{Q}[G]$ corresponding to the extremely strong Shoda pairs of $G$ in $\mathcal{S}_{N}$ is added to $SumDim$. In view of Lemma \ref{l3}, the process stops when either $SumDim=|G|$ or when all the normal subgroups of $G$ are exhausted. The normal subgroups $N$ of $G$ are selected in decreasing order i.e., if the normal subgroup $N_{1}$ is chosen before the normal subgroup $N_{2}$, then $|N_{1}|\geq|N_{2}|$. This has been done keeping in view the ease of computation of $\mathcal{S}_{N}$, if $G/N$ has small order. This algorithm enables us to write the function \texttt{ExtStrongShodaPairs(G);} in \texttt{GAP} language.

 \subsection{Strong Shoda Pairs}

 We next describe Algorithm \ref{A2} to compute the set $StSP$, which is a complete irredundant set of strong Shoda pairs of a given finite group $G$.\\

\begin{algorithm}[H]\label{A2}
\begin{scriptsize}
\KwData{A finite group $G$.}
$StSP$:= A complete irredundant set of extremely strong Shoda pairs of $G$;\\
$SumDim$:=the sum of $\mathbb{Q}$-dimensions of simple components of $\mathbb{Q}[G]$ corresponding to the primitive central idempotents realized by the extremely strong Shoda pairs of $G$;

\eIf{$SumDim$=$|G|$}{return $StSP$;}{

$PCIs$:= Primitive central idempotents of $\mathbb{Q}[G]$ realized by strong Shoda pairs in $StSP$;\\
$C$:= Conjugacy classes of those subgroups of $G$ which are not normal;\\
\While{ $SumDim \neq |G|$}{

\For{$c$ in $C$}{$K:=$Representative$(c);$ $N:=\operatorname{core}_{G}(K)$;

\If{$Centre(G/N)$ is cyclic}{ $H$:=a subgroup of $G$ such that $(H,K)$ is a strong Shoda pair and $H \ntrianglelefteq G$;\\ $e:=e(G,H,K);$\\
\If{$e$ is not in the list $PCIs$}{Add $e$ to the list $PCIs$;\\  Add $[H,K]$ to the list $StSP$;\\ $dim:=\operatorname{dim}(H,K);$ $SumDim:=SumDim+dim;$}}
 }}}

\KwResult{$StSP$}
\caption{Strong Shoda pairs of $G$}

\end{scriptsize}
\end{algorithm}
\vspace{0.5cm}
\para
 Initially, $StSP$ is the list $ESSP$ of extremely strong Shoda pairs of $G$ obtained using Algorithm \ref{A1}. Also, $SumDim$ is set to be the the sum of $\mathbb{Q}$-dimensions of simple components of $\mathbb{Q}[G]$ corresponding to the primitive central idempotents realized by extremely strong Shoda pairs of $G$. In case $SumDim=|G|$, by Corollaries \ref{c1} and \ref{c2}, $StSP$ is a complete irredundant set of strong Shoda pairs of $G$ and the algorithm terminates. Otherwise, to find the remaining strong Shoda pairs of~$G$, we make use of the algorithm provided by  Olivieri and del R{\'{\i}}o \cite{OliA} with desired modifications. For a strong Shoda pair $(H,K)$ of $G$, we use the fact that $G/\operatorname{core}_{G}(K)$ must be cyclic (Lemma \ref{l2}). Moreover, if $(H,K)$ realizes a primitive central idempotent of $\mathbb{Q}[G]$ different from the one realized by an extremely strong Shoda pair of $G$, then none of $H$ or $K$ is normal in $G$. This algorithm allows us to write the function \texttt{StShodaPairs(G);} in \texttt{GAP} language.

\subsection{Primitive Central Idempotents}
The algorithm to compute the primitive central idempotents of $\mathbb{Q}[G]$ realized by extremely strong Shoda pairs of $G$ is similar to Algorithm 1. The only difference is that at any stage of the computation, instead of collecting the elements of $\mathcal{S}_{N}$, one collects the primitive central idempotents realized by them. Using this algorithm, we write the function \texttt{PrimitiveCentralIdempotentsByExtSSP(QG);} in \texttt{GAP} language which computes the set of primitive central idempotents realized by extremely strong Shoda pairs of $G$. To compute the primitive central idempotent $e(G,H,K)$ of $\mathbb{Q}[G]$ realized by the strong Shoda pair $(H,K)$ of $G$, we use the function \texttt{Idempotent\_eGsum(QG,H,K);} currently available in \texttt{Wedderga}.

\para Similarly, the algorithm to compute the primitive central idempotents of $\mathbb{Q}[G]$ realized by strong Shoda pairs of $G$ is obtained by a slight modification of Algorithm~\ref{A2} and the corresponding function \texttt{PrimitiveCentralIdempotentsByStSP(QG);} is also obtained.

\subsection{Normally Monomial Groups}

The algorithm to check whether a finite group $G$ is normally monomial or not is obtained by replacing the result $ESSP$ of Algorithm \ref{A1} with $SumDim$. In view of Corollary \ref{c2}, $G$ is normally monomial if, and only if, $SumDim$=$|G|$. This algorithm enables us to write the function \texttt{IsNormallyMonomial(G);} in \texttt{GAP} language.
 \para Using the function \texttt{IsNormallyMonomial(G);} we have found by a computer search that $98.84\%$ of the monomial groups of order up to $500$ are normally monomial. Also, $97.88\%$ of all the finite groups of order up to $500$ are normally monomial. An exhaustive computer search also yields that among the groups of odd order up to $2000$, the only groups which are not normally monomial are:

\begin{quote}
   \texttt{SmallGroup(375,2); SmallGroup(1029,12); SmallGroup(1053,51); SmallGroup(1125,3); SmallGroup(1125,7); SmallGroup(1215,68); SmallGroup(1875,18); SmallGroup(1875,19);}
\end{quote}

\noindent It may be pointed out that all the groups in the above list, except the second and the third, are non monomial.\\

\section{Runtime Comparison}

  \para We now present an experimental runtime comparison between the  following two sets of functions for different samples of groups:
   \begin{enumerate}
    \item  \texttt{StrongShodaPairs(G);} with \texttt{StShodaPairs(G);}
    \item  \texttt{PrimitiveCentralIdempotentsByStrongSP(QG);} with \\ \texttt{PrimitiveCentralIdempotentsByStSP(QG);}
   \end{enumerate}

\para For a given sample of groups, let $t(n)$ be the average of the runtimes, taken in milliseconds, for the groups in $S$ of order $n$ for $n\geq 1$. If the sample contains no group of order $n$, then set $t(n)=0$. Define $T(n)= \sum_{i=1}^{n}t(i),~n\geq 1$.

 \para We now describe the first sample $S$ which consists of $31272$ groups of order up to $2000$. For $1 \leq n \leq 2000$, $n\neq 1024$, if the number of non isomorphic groups of order $n$ is less than $200$, then $S$ contains all the groups of order $n$. Otherwise, we include in the sample $S$, at least $100$ groups of order $n$, which are evenly spread in the \texttt{GAP} library of small groups. The groups of order 1024 are excluded because of their non availability in \texttt{GAP} library. For this sample, the graph of~$n$ versus $T(n)$ for the comparison of the functions \texttt{StrongShodaPairs(G);} and \texttt{StShodaPairs(G);} is presented in Fig.1. Also, Fig.2 presents the runtime comparison of the function \texttt{PrimitiveCentralIdempotentsByStrongSP(QG);} with the function \texttt{PrimitiveCentralIdempotentsByStSP(QG);}

\begin{figure}[H]\label{f1}\centering
\includegraphics[scale=.65]{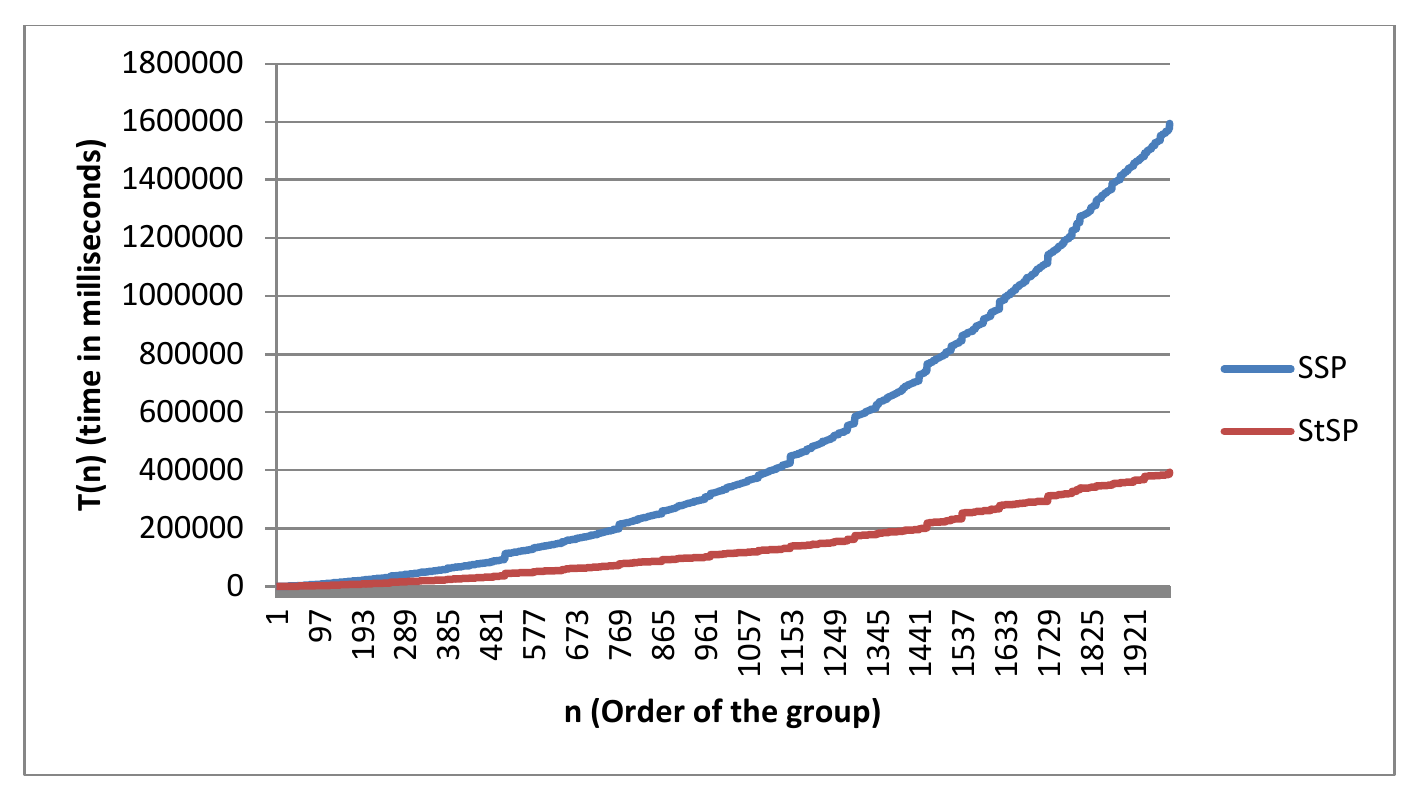}
\renewcommand
{\figurename}{Fig.}
\caption{Strong Shoda pairs~(Sample S)}
\end{figure}

\begin{figure}[H]\label{f2}
\centering
\includegraphics[scale=.65]{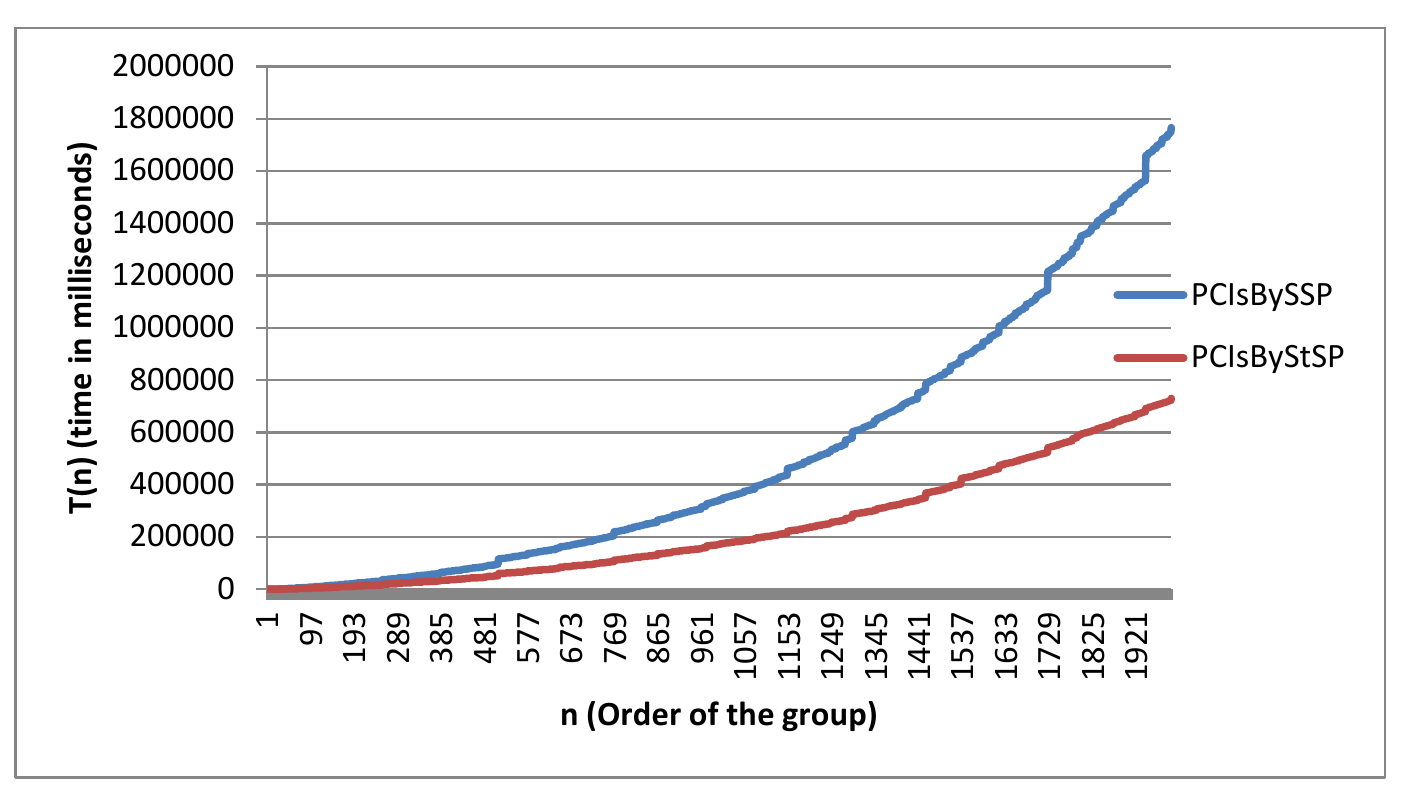}
\renewcommand
{\figurename}{Fig.}
\caption{Primitive Central Idempotents~(Sample S)}
\end{figure}

We next compare the runtimes of \texttt{StrongShodaPairs(G);} and \texttt{StShodaPairs(G);} for a sample of solvable and that of non solvable groups. The sample $S_{1}$ of solvable groups consists of all the groups of odd order up to 2000 and the sample $S_{2}$ consists of all the non solvable groups of order up to 2000. The graph of $n$ versus $T(n)$ for these samples are presented in Figs.3 and 4 respectively.

\begin{figure}[H]\label{f3}\centering
\includegraphics[scale=.65]{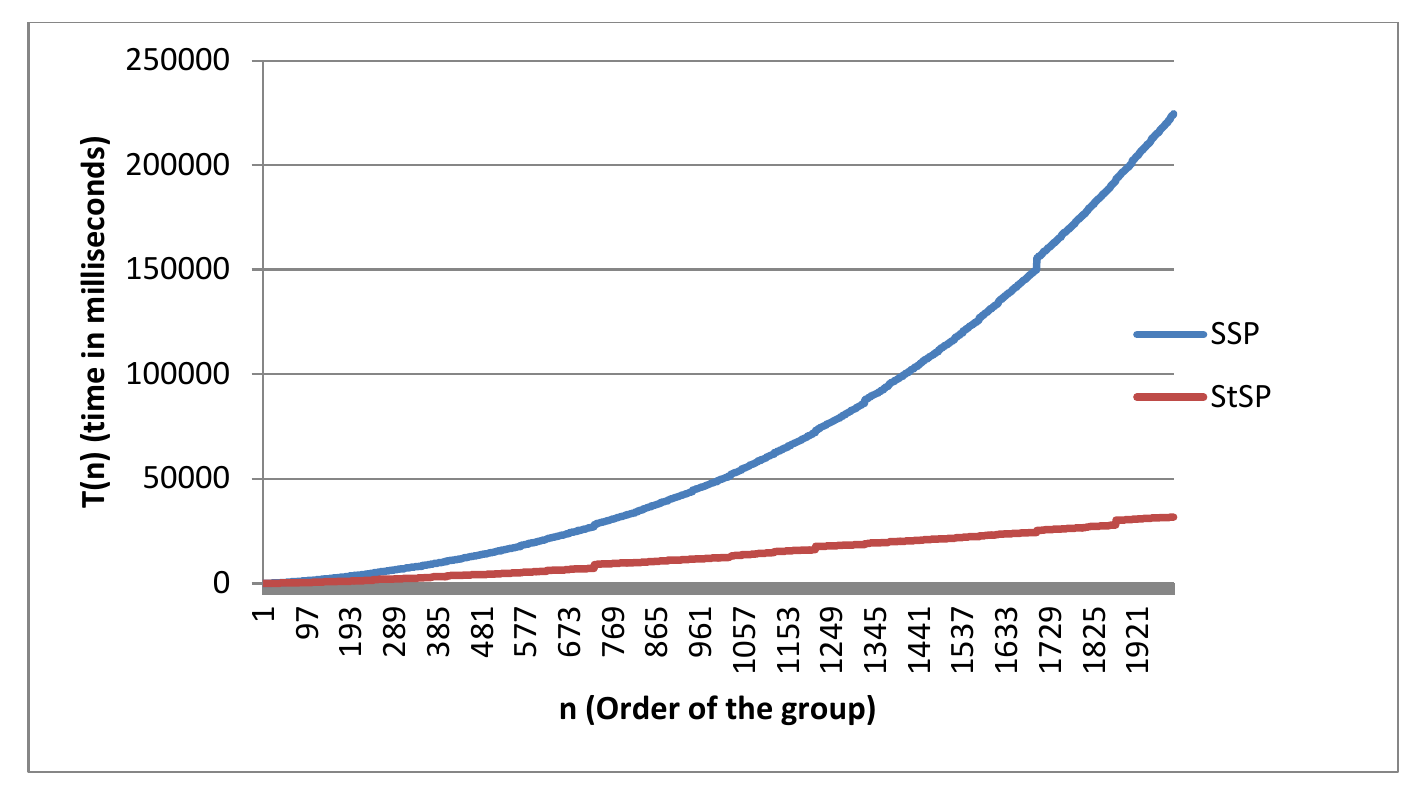}
\renewcommand
{\figurename}{Fig.}
\caption{Strong Shoda pairs~(Sample S1)}
\end{figure}

\begin{figure}[H]\label{f4}
\centering
\includegraphics[scale=.65]{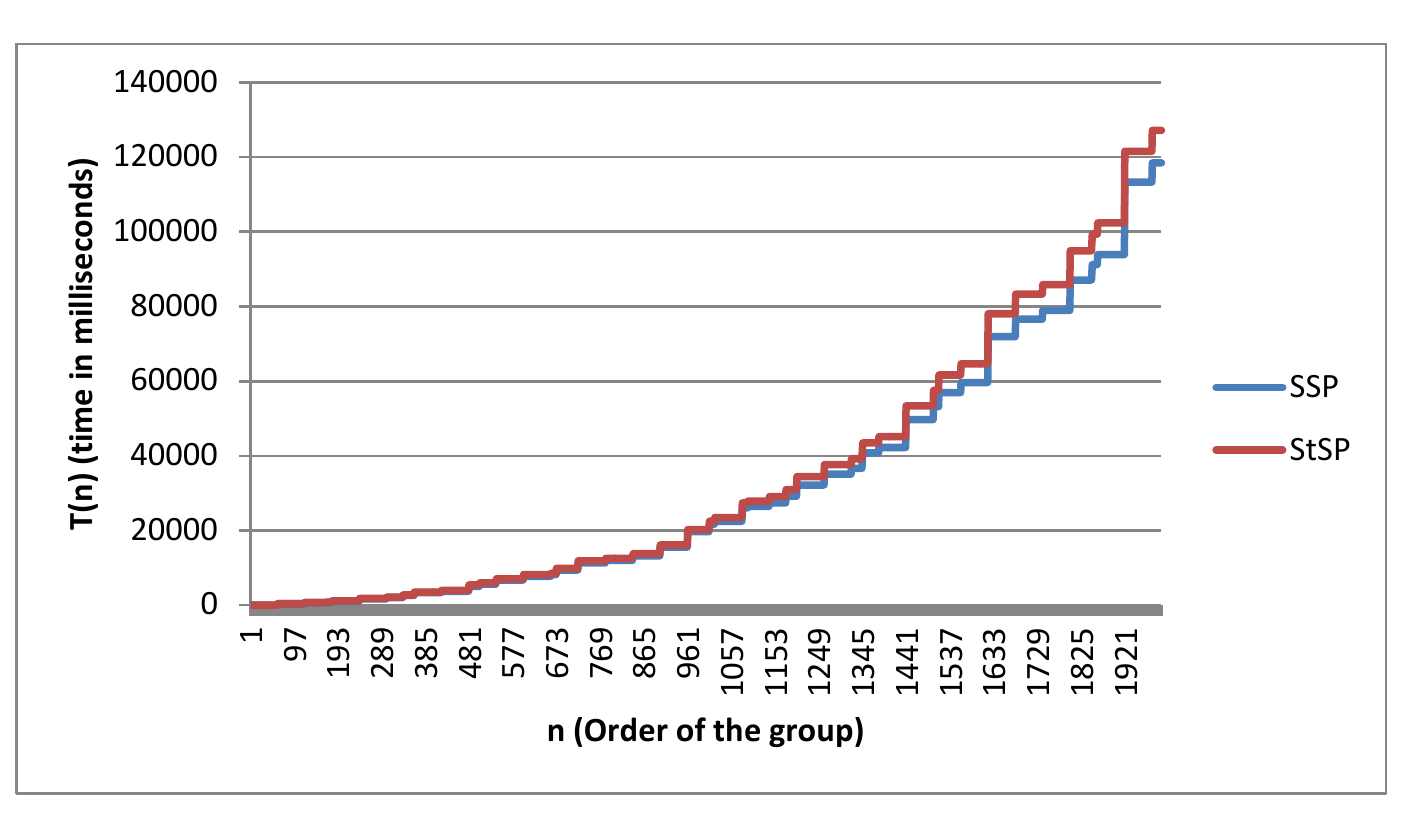}
\renewcommand
{\figurename}{Fig.}
\caption{Strong Shoda pairs~(Sample S2)}
\end{figure}

\para In Figs.1-4, \texttt{SSP} and \texttt{StSP} are the curves for the functions \texttt{StrongShodaPairs(G);} and \texttt{StShodaPairs(G);} respectively and \texttt{PCIsBySSP} and \texttt{PCIsByStSP} represent the curves for the functions \texttt{PrimitiveCentralIdempotentsByStrongSP(QG);} and \texttt{PrimitiveCentralIdempotentsByStSP(QG);} respectively. These experiments have been performed on the computer with Intel Core i7-4770 CPU @ 3.40GHz Dual Core, 4GB RAM.\\

 The overall improvement in the performance of \texttt{StShodaPairs(G);} in comparison to \texttt{StrongShodaPairs(G);} is mainly due to following differences in their respective algorithms:

\begin{multicols}{2}

\begin{center}\texttt{StShodaPairs(G);}

\columnbreak

\texttt{StrongShodaPairs(G);}\end{center}

\end{multicols}

\begin{multicols}{2}
\begin{description}
\item[1.] \hspace{.45cm}Begins by computing all the normal subgroups of $G$.  The conjugacy classes of subgroups of $G$ are computed only if $G$ is not normally monomial.
\end{description}
\columnbreak
\begin{description}
\item[1.] \hspace{.45cm}Always begins by computing all conjugacy classes of subgroups of~$G$. It may be pointed out that generating the full subgroup lattice of $G$ restricts the efficiency when $G$ has large order.
    \end{description}

\end{multicols}

\begin{multicols}{2}
\begin{description}
  \item[2.]\hspace{.29cm}  Firstly, the extremely strong Shoda pairs of $G$ are computed. If $G$ is not normally monomial, then the remaining strong Shoda pairs of are found by the search algorithm of \texttt{StrongShodaPairs(G);}, with slight modifications.
\end{description}

\columnbreak
\begin{description}
  \item[2.]\hspace{.2cm} There is no distinction between the computation of extremely strong Shoda pairs and that of strong Shoda pairs of $G$.
\end{description}
\end{multicols}

\begin{multicols}{2}
\begin{description}
  \item[3.]\hspace{.13cm} Extremely strong Shoda pairs of $G$ are computed using Theorem\ref{t2}, which ensures that each time a new extremely strong Shoda pair is constructed, it is necessarily inequivalent to any of the extremely strong Shoda pair already obtained.
\end{description}

\columnbreak
\begin{description}
  \item[3.] \hspace{.45cm}When a new strong Shoda pair of $G$ is discovered, it is not necessarily inequivalent to the ones already discovered. Each time a new strong Shoda pair is found, the algorithm computes the corresponding primitive central idempotent of $\mathbb{Q}[G]$ and checks its equivalence.
\end{description}

\end{multicols}

\vspace{0.5cm}
The above differences also result in the improved performance of the function \texttt{PrimitiveCentralIdempotentsByStSP(QG);} in comparison to that of the function  \texttt{PrimitiveCentralIdempotentsByStrongSP(QG);} which is currently available in \texttt{Wedderga}.

\vspace{0.7cm}
\noindent\textbf{Acknowledgements}\\

\noindent The authors are grateful to the anonymous referees for their valuable comments and suggestions which have helped to write the paper in the present form.

\bibliographystyle{amsplain}
\bibliography{ReferencesBM3}
\end{document}